\documentclass[twoside,11pt,reqno]{amsart}

\usepackage{amsmath,amssymb,amscd,mathrsfs,amscd, todonotes,appendix}
\usepackage{graphics,verbatim}
\usepackage{todonotes}
\usepackage{enumitem}
\usepackage{hyperref}

\makeatletter
\@namedef{subjclassname@2020}{\textup{2020} Mathematics Subject Classification}
\makeatother

\newcommand{\gl}{\mathfrak{g}}
\newcommand{\opH}{\operatorname{H}}
\newcommand{\ul}{\mathfrak{u}}
\newcommand{\ind}{\operatorname{ind}}
\newcommand{\gfp}{G(\mathbb{F}_p)}
\newcommand{\gfq}{G(\mathbb{F}_{q})}

\newcommand{\la}{\lambda}
\newcommand{\si}{\sigma}

\newcommand{\al}{\alpha}
\newcommand{\be}{\beta}
\newcommand{\Ext}{\operatorname{Ext}}
\newcommand{\Hom}{\operatorname{Hom}}

\newcommand{\fp}{{\mathbb F}_p}
\newcommand{\fq}{{\mathbb F}_q}
\newcommand{\res}{\rm{res}}

\makeatletter
\newcommand{\leqnomode}{\tagsleft@true}
\newcommand{\reqnomode}{\tagsleft@false}
\makeatother

\makeatletter\let\c@fact\c@theorem\makeatother

\makeatletter\let\c@note\c@theorem\makeatother

\makeatletter\let\c@lemma\c@theorem\makeatother

\makeatletter\let\c@lemma\c@theorem\makeatother

\newtheorem{quest}{Question}[subsection]
\makeatletter\let\c@alg\c@theorem\makeatother

\makeatletter\let\c@prop\c@theorem\makeatother

\makeatletter\let\c@conj\c@theorem\makeatother

\makeatletter\let\c@cor\c@theorem\makeatother

\makeatletter\let\c@defn\c@theorem\makeatother

\theoremstyle{definition}

\newtheorem{remark}{Remark}[subsection]
\makeatletter\let\c@remark\c@theorem\makeatother

\makeatletter\let\c@example\c@theorem\makeatother
\numberwithin{equation}{subsection}

\usepackage[capitalise]{cleveref}
\crefname{theorem}{Theorem}{Theorems}
\crefname{fact}{Fact}{Facts}
\crefname{note}{Note}{Notes}
\crefname{lemma}{Lemma}{Lemmas}
\crefname{alg}{Algorithm}{Algorithms}
\crefname{remark}{Remark}{Remarks}
\crefname{example}{Example}{Examples}
\crefname{prop}{Proposition}{Propositions}
\crefname{conj}{Conjecture}{Conjectures}
\crefname{cor}{Corollary}{Corollaries}
\crefname{defn}{Definition}{Definitions}
\crefname{equation}{\!\!}{\!\!} 

\newcounter{listequation}


\begin{document}

\title{Cohomology of algebraic groups, Lie algebras, and finite groups of Lie type}

\dedicatory{In memory of Georgia M. Benkart and Brian J. Parshall}

\begin{abstract} Let $G$ be a reductive algebraic group over a field of prime characteristic.  One can associate to $G$ (or subgroups thereof) its Lie algebra, its Frobenius kernels, and the finite Chevalley group of points over a finite field.  The representation theories of these structures are highly interconnected.  This expository article will focus specifically on the cohomology theories of these structures and the relationships between them with the aim of highlighting a few key developments over the past 20 years and related open questions. 
\end{abstract}

\author{\sc Christopher P. Bendel}
\address
{Department of Mathematics, Statistics and Computer Science\\
University of
Wisconsin-Stout \\
Menomonie\\ WI~54751, USA}
\thanks{Research of the author was supported in part by Simons Foundation Collaboration Grant 317062}
\email{bendelc@uwstout.edu}

\maketitle

\section{Introduction}
		
\subsection{}\label{S:intro} Let $k$ be an algebraically closed field of characteristic $p > 0$ and $G$ be a connected, reductive algebraic group scheme over $k$.  For example, $G$ may be the general linear group $GL_n$ or the (semi-simple) special linear group $SL_n$ or one of the other classical matrix groups.    Associated to $G$ are several standard subgroup schemes: a maximal torus $T$, a Borel subgroup $B$ containing $T$, and the unipotent radical $U$ of $B$.  For example, in the case of $GL_n$, $T$ consists of the diagonal matrices, $B$ might be taken to be the (non-strictly) lower triangular matrices, and then $U$ is the strictly lower triangular matrices.   

There are a number of other algebraic structures related to $G$: Frobenius kernels, Lie algebras, and finite groups of Lie type.   Let $F: G \to G$ denote the Frobenius map (given in the matrix example by raising each matrix entry to the $p$th power).   The first Frobenius kernel of $G$ is the scheme-theoretic kernel of $F$, denoted $G_1$. More generally, for an integer $r \geq 1$, the $r$th Frobenius kernel $G_r$ is the kernel of the $r$th iterate $F^r$ of $F$.  One similarly has Frobenius kernels for the aforementioned subgroup schemes $B$ and $U$.    

The first Frobenius kernel $G_1$ is closely related to the Lie algebra $\gl$ of $G$.   In this modular setting, $\gl$ admits the additional structure of a $p$-restricted Lie algebra.   In the matrix example, the defining $[p]$-power mapping is given by $p$th power of matrices.   Similarly, one has ($p$-restricted) Lie algebras associated to $B$ and $U$.  The representations of $G_1$ are equivalent to the restricted representations of $\gl$.  Formally, there is an isomorphism of Hopf algebras $k[G_1]^* \cong u(\gl)$, where $k[G_1]^*$ denotes the $k$-linear dual of the coordinate algebra of $G_1$ and $u(\gl)$ denotes the restricted enveloping algebra of $\gl$.  Geometrically, one may think of the Lie algebra $\gl$ as a ``tangent space'' approximation to $G$.   Translating that to Frobenius kernels, as $r$ increases, the representation theory of $G_r$ provides a better approximation to that of $G$.  See Section \ref{S:mod} for an example of this behavior.  

Let $\fp$ denote the finite field of $p$ elements and, more generally, $\fq$ the field with $q := p^r$ elements.  The fixed  points of the Frobenius map (or its iterates) give rise to finite groups of Lie type: $\gfp := G^F$ and $\gfq := G^{F^r}$.  In our matrix example, these are the finite matrix groups with entries in the given finite field.  One similarly has finite groups associated to $B$ and $U$.  In addition, one has ``twisted'' groups of Lie type that arise by composing the Frobenius map with an automorphism of $G$ induced from a non-trivial automorphism of the associated Dynkin diagram.  See for example \cite[1.5]{Hum}.  Many of the results presented in this paper hold also for these twisted groups. For convenience, the discussion here will be restricted to the ``untwisted'' case.  

Any rational $G$-module may be restricted to a module for any of the above objects, and the representation theories of these structures are well-known to be highly interconnected.  For a thorough discussion, the reader is directed to the excellent books by Jantzen \cite{Jan} and Humphreys \cite{Hum}, noting that the notation used here generally follows that of \cite{Jan}.  The goal of this paper is not to reproduce the material found therein, rather it is an incomplete attempt to highlight the current status of a few fundamental cohomological questions where progress has been made since those books were published. It is also hoped that this work complements and updates the article of Nakano \cite{Nak}.  As done in \cite{Nak}, this work aims to include examples demonstrating how the cohomology of these various structures are interrelated.

Sections \ref{S:trivial} and \ref{S:nontriv} consider the cohomology of Frobenius kernels, beginning with trivial coefficients in Section \ref{S:trivial}.   Computations of of $\opH^i(G_r,k)$, $\opH^i(B_r,k)$, and $\opH^i(U_r,k)$ in particular degrees will be considered as well as the global ring structure of $\opH^{\bullet}(G_r,k) := \bigoplus_{i \geq 0}\opH^i(G_r,k)$ (or for $B$ or $U$).  In Section \ref{S:nontriv}, non-trivial coefficients are considered, with the primary focus being on the computation of the groups $\opH^i(G_r,H^0(\la))$ for a standard induced module $H^0(\la)$, along with the related computation of $\opH^i(B_r,\la)$ for the one-dimensional $B$-module with weight $\la$. See Section \ref{S:mod} for details on the module notation.  This also leads to a discussion of the ordinary Lie algebra cohomology of $\ul$, i.e., $\opH^i(\ul,k)$.    

The discussion for Frobenius kernels leads to algebraic group computations, including a discussion of $B$-cohomology $\opH^i(B,\la)$ in Section \ref{S:alg} and of $G$-cohomology in Section \ref{S:rat} with a new look at generic cohomology and the connections between $G$-cohomology and $\gfq$-cohomology.   Cline, Parshall, Scott, and van der Kallen \cite{CPSvdK} showed that the cohomology of the finite group $\gfq$ stabilizes as $r$ increases and may be identified with certain $G$-cohomology.   Section \ref{S:rat} will discuss recent improvements to that result and further connections between $\gfq$- and $G$-cohomology for simple modules.  These recent results are highly dependent on the use of the induction functor $\ind_{\gfq}^G(-)$, a tool that was introduced into the theory by Bendel, Nakano, and Pillen and has proved to serve as a way to conceptually unify a number of finite group cohomology questions, as will be seen in applications in Sections \ref{S:rat}, \ref{S:van}, \ref{S:small}, and \ref{S:bound}.  

Section \ref{S:van} will continue to discuss $\gfq$-cohomology, bringing in ideas from Sections \ref{S:nontriv} and \ref{S:rat}.  Dating to work of Quillen \cite{Q}, it has been observed that $\opH^i(\gfq,k)$ is zero in ``small'' positive degrees.   This section will discuss results giving the least positive degree with non-trivial cohomology. Low degree $\gfq$-cohomology of simple modules is considered in Section \ref{S:small}.  Lastly, Section \ref{S:bound} will consider, in multiple contexts, the problem of determining the maximum dimension of a cohomology group.  Ideas from previous sections will again play a role in the results presented here.

This paper will primarily consider cohomology, which can of course be viewed from the perspective of extension groups. E.g., for a rational $G$-module $M$, $\opH^i(G,M) \cong \Ext_{G}^i(k,M)$.   In some cases, explicit mention of related Ext-questions will be made, with other generalizations left to the interested reader.    With this commentary on extensions, the reader may note that absent from the above list of topics is the Lusztig Conjecture on the character of a simple $G$-module; a topic with strong extension/cohomological ties.   There have been many exciting developments on that question in recent years, and a thorough discussion thereof is beyond the scope of this article. The interested reader is directed to the recent work of Riche and Williamson \cite{RW} (and references therein).

\subsection{Notation}\label{S:not} From now on, except as otherwise noted, $G$ will denote a simple, simply connected algebraic group which is defined and split over $\fp$.  As above, $T$ denotes a maximal torus.   Associated to $G$ is an irreducible root system $\Phi$ with positive roots $\Phi^+$ and simple roots $S$.   $B$ will be a Borel subgroup containing $T$ corresponding to the negative roots, with $U$ the unipotent radical so that $B = T\ltimes U$.

Associated to the root system $\Phi$, the following notation will be used:
\begin{itemize}
\item $\alpha_0$: the maximal short root.

\item $\mathbb E$: the Euclidean space spanned by $\Phi$ with inner product $\langle\,,\,\rangle$ normalized so that $\langle\alpha,\alpha\rangle=2$ for $\alpha \in \Phi$ any short root.

\item $\alpha^\vee=2\alpha/\langle\alpha,\alpha\rangle$: the coroot of $\alpha\in \Phi$.

\item $\rho$: the Weyl weight defined by $\rho=\frac{1}{2}\sum_{\alpha\in\Phi^+}\alpha$.

\item  $h$: the Coxeter number of $\Phi$, given by $h=\langle\rho,\alpha_0^{\vee} \rangle+1$. Explicit values are given here:

\vskip.2cm
\begin{center}
\begin{tabular}{|c|c|c|c|c|c|c|c|c|c|}
\hline
$\Phi$ & $A_n$ & $B_n$ & $C_n$ & $D_n$ & $E_6$ & $E_7$ & $E_8$ & $F_4$ & $G_2$\\
\hline
$h$ & $n + 1$ & $2n$ & $2n$ & $2n - 2$ & 12 & 18 & 30 & 12 & 6\\
\hline
\end{tabular}
\end{center}
\vskip.2cm

\item $W$: the Weyl group of $\Phi$. 

 
\item $\ell : W \to\mathbb N$: the usual length function on $W$.

\item $X=\mathbb Z \varpi_1\oplus\cdots\oplus{\mathbb Z}\varpi_n$: the weight lattice, where the fundamental dominant weights $\varpi_i\in{\mathbb E}$ are defined by $\langle\varpi_i,\alpha_j^\vee\rangle=\delta_{ij}$, $1\leq i,j\leq n$.

\item $X^+={\mathbb N}\varpi_1+\cdots+{\mathbb N}\varpi_n$: the cone of dominant weights.

\item  $X_{r}=\{\lambda\in X^+: 0\leq \langle\lambda,\alpha^\vee\rangle<p^{r},\,\,\forall\, \alpha\in S\}$: the set of $p^{r}$-restricted dominant weights. 



\item $\lambda^\star : = -w_0\lambda$: where $w_0$ is the longest word in the Weyl group $W$ and $\lambda\in X$. 

\item $\leq$: the usual partial ordering on $X$, where $\la \leq \mu$ iff $\la - \mu$ is a non-negative sum of simple roots.

\item The dot action $\cdot$ on $X$: for $w \in W$ and $\la \in X$, $w\cdot \la = w(\la + \rho) - \rho$.
\end{itemize}

\subsection{Modules}\label{S:mod} Let $H$ be an arbitrary affine algebraic group scheme over $k$ and $M$ be a rational $H$-module. $M^*$ will denote the ordinary $k$-linear dual of $M$.  Let $M^{(r)}$ denote the ``twisted'' representation obtained by composing the underlying representation of $H$ with $F^r$ (the $r$th iterate of the Frobenius morphism on $H$).   If one has an isomorphism $M \cong N^{(r)}$ of rational $H$-modules, then one sometimes writes $M^{(-r)}$ for the module $N$.  

Continuing with our notation as above, consider the Borel subgroup $B$ and a weight $\la \in X$.   By definition, such a weight $\la$ defines a one-dimensional $T$-module.  Letting $U$ act trivially gives rise to a one-dimensional $B$-module that will also (abusively) be denoted $\la$.  Induction gives rise to one of our standard $G$-modules: set $H^0(\la) := \ind_B^G\la$.  It is well-known that $H^0(\la) \neq 0$ iff $\la \in X^+$.   When $\la$ is dominant, $H^0(\la)$ has highest weight $\la$ and simple socle, denoted $L(\la)$.    The set of such $L(\la)$ over $X^+$ gives a complete set of finite-dimensional simple $G$-modules.   The induced module $H^0(\la)$ is dual to a Weyl module. Let $V(\la)$ denote the Weyl module with highest weight $\la$.  One has $V(\lambda)\cong H^0(\lambda^\star)^*$ and the character of $V(\la)$ or $H^0(\la)$ is given by Weyl's character formula (cf. \cite[II.5.10,5.11]{Jan}).

Given a rational $G$-module $M$, one may restrict to get a module over $B$, $U$, $G_r$, $B_r$, $U_r$, and $\gfq$ (as well as subgroups thereof).  In particular, the finite-dimensional simple modules for $G_r$ and $\gfq$ arise in such a way.   Given $\la \in X_r$, the restriction of $L(\la)$ to $G_r$ or $\gfq$ remains simple and the set of all such $L(\la)$ gives a complete set of finite-dimensional simple modules.   This is perhaps the most fundamental example of the connections between the representation theories of $G_r$ and $\gfq$, and of both groups to $G$.  This also illustrates the conceptual idea noted in Section \ref{S:intro} that $G_r$-representation theory approximates that of $G$: as $r$ increases the set of $G_r$-simples approaches the set of $G$-simples.

\subsection{Spectral sequences} Various spectral sequences have proven to be key tools in the study of cohomology.  Most fundamental is the Lyndon-Hochschild-Serre spectral sequence associated to a normal subgroup of a group. See for example \cite[Prop. I.6.6]{Jan}.  

A second spectral sequence gives an immediate connection between the (ordinary) cohomology of Lie algebras and the cohomology of a first Frobenius kernel (cf. \cite[Rem. I.9.20]{Jan}). Let $H$ be an arbitrary affine algebraic group scheme over $k$ and assume $p \neq 2$, with $\mathfrak{h}$ denoting the Lie algebra of $H$.  Let $M$ be an $H_1$-module, then there exists a first-quadrant spectral sequence with $E_2$-term as follows

\begin{equation}\label{E:AJss}
E_2^{2i,j} = S^i(\mathfrak{h}^*)^{(1)}\otimes \opH^j(\mathfrak{h},k) \Rightarrow \opH^{2i + j}(H_1,M).
\end{equation}

A third relevant spectral sequence relates $G$-cohomology to $B$-cohomology via the induction functor.  For a $B$-module $M$ such that $R^m\ind_B^G(M) = 0$ for all $m > 0$, there exists a first quadrant spectral sequence (cf. \cite[II.12.2]{Jan})
\begin{equation}\label{E:indss}
E_2^{i,j} = R^i\ind_B^G\left(\opH^j(B_r,M)^{(-r)}\right)\Rightarrow \opH^{i+j}\left(G_r,\ind_{B}^{G}(M)\right)^{(-r)}.
\end{equation}
This may be applied for example to $M = \la$ for a dominant weight $\la$.

\subsection{Acknowledgements} I would like to thank Daniel Nakano and Cornelius Pillen for their many years of friendship and mathematical collaboration that led to a number of the results discussed in this paper.  I also thank the organizers for the opportunity to speak at the Southeastern Lie Theory Workshop X held in June 2018 at the University of Georgia.  This article stemmed from the material given in a series of of lectures presented at that meeting.  

The year 2022 saw the tragic passing of two stalwarts in the field: Georgia Benkart and Brian Parshall.   Both were in attendance at that SE Lie Theory Workshop, and Georgia was my ``partner'' in giving a second minicourse at the event. I would like to acknowledge their impact on my mathematical career.  I received my first speaking invitation from Georgia, and we interacted regularly over the years, particularly as she tried to promote representation theory in Wisconsin and the surrounding region.   I also met Brian early in my career, when he visited Northwestern University.  I eventually had the pleasure of collaborating with him on multiple projects.  I remain ever grateful for a wonderful semester-long sabbatical I spent at the University of Virginia; time that bore mathematical fruit for many years thereafter.


\section{Cohomology of Frobenius Kernels with Trivial Coefficients}\label{S:trivial}

\subsection{} Consider the Frobenius kernel $G_r$.   Of interest here is not only the cohomology group $\opH^i(G_r,k)$ for a given degree but also the cohomology ring 
$$
\opH^{\bullet}(G_r,k) := \bigoplus_{i = 0}^{\infty}\opH^i(G_r,k).
$$
This is a graded commutative ring and, by work of Friedlander and Suslin \cite{FS}, known to be finitely-generated. While these rings have been described geometrically in a sense through the use of cohomological varieties\footnote{Another topic with an extensive literature that will not be discussed here.} (cf. \cite{SFB1, SFB2}), explicit ring computations remain minimal.

\subsection{} In the $r = 1$ case, for sufficiently large $p$, identification of $\opH^{\bullet}(G_1,k)$ preceded the general finite-generation result of Friedlander-Suslin.    It was shown by Friedlander and Parshall \cite{FP86} and also by Andersen and Jantzen \cite{AJ} that, for $p > h$, there is an algebra isomorphism
$$
\opH^{\bullet}(G_1,k)^{(-1)} \cong k[{\mathcal N}],
$$
where ${\mathcal N} \subseteq \gl$ is the set of nilpotent elements in $\gl$.   A key idea in the proof was the use of induction from $B$ to $G$.  

\subsection{The Induction Question}\label{S:ind} Recall the induction spectral sequence relating $G_r$-cohomology to $B_r$-cohomology given in \eqref{E:indss}. 
Optimally, one would have an affirmative answer to the following key question (for a $B$-module $M$ with $R^i\ind_B^G(M) = 0$ for $i > 0$):

\begin{quest}\label{Q:ind} Is $R^i\ind_{B}^{G}\left(\opH^j(B_r,M)^{(-r)}\right) = 0$ for $i > 0$?  
\end{quest}

If yes, the spectral sequence would collapse to give 
$$
\opH^j\left(G_r,\ind_{B}^{G}(M)\right)^{(-r)} \cong \ind_{B}^{G}\left(\opH^j(B_r,M)^{(-r)}\right).
$$
For both answering the question and then making a final computation, it is often useful to use the fact that $\opH^i(B_r,M) \cong \opH^i(U_r,M)^{T_r}$.

Indeed, for $M = k$, this strategy works to give the above computation of $G_1$-cohomology. The condition $p > h$ arises when considering $T_1$-fixed points.   See Section \ref{S:B1} below for more on this.  This approach will be applied at multiple points below, and Question \ref{Q:ind} remains a highly relevant question in its own right, in addition to its potential applications.

\subsection{Small primes} For small primes, examples in \cite{AJ} show that the answer is necessarily different, with the nilpotent cone (or nullcone)  ${\mathcal N}$ needing to be replaced by a smaller variety.   Observe that ${\mathcal N}$ may be identified with $G\cdot \ul$, the $G$-orbit of the Lie algebra of $U$.  When $2 < p \leq h$ (and some further mild restrictions on $p$), work of Bendel, Nakano, Parshall, and Pillen \cite{BNPP}\footnote{The primary result of that work was an identification of the complex cohomology ring of what is often known as Lusztig's ``small'' quantum group for $\gl$ at a root of unity.} provides an answer in some cases and a conjectural answer more generally in terms of certain $G$-orbits.  

Assume that $p \geq 3$ and that $p$ is a {\em good} prime for the root system.  More precisely, the latter condition requires that $p \geq 5$ in types $E_6$, $E_7$, $F_4$, and $G_2$ and $p \geq 7$ in type $E_8$.   In type $A_n$, assume also that $p$ does not divide $n + 1$ ($p$ is said to be {\em very good}).  Consider the subset
$$
\Phi_0 := \{\al \in \Phi ~ | ~ \langle \rho, \al\rangle \equiv 0 \mod p\} \subseteq \Phi.
$$
When $p \geq h$, $\Phi_0$ is empty, whereas, when $p < h$, $\Phi_0$ is a closed nonempty subsystem of $\Phi$.  There exists a subset $J \subseteq S$ of simple roots and an element $w \in W$ such that $w(\Phi_0) = \Phi_J$ (the root subsystem generated by $J$) and $w(\Phi_0^+) = \Phi_J^+$.  Associated to $J$ is a standard parabolic subgroup $B \subseteq P_J = L_J\ltimes U_J \subseteq G$ with Levi subgroup $L_J$ and unipotent radical $U_J$.  Let $\ul_J$ denote the Lie algebra of $U_J$.  Note that, if $J$ was the empty set, we would have $P_J = B$, $L_J = T$, and $U_J = U$.  

One can replace the spectral sequence of \eqref{E:indss} with one involving the parabolic subgroup $P_J$ and the question of computing $T_1$ fixed points with understanding the following homomorphism group:
$$
\Hom_{(L_J)_1}\left(k,\ind_B^{P_J}w\cdot 0\otimes H^j(\ul_J,k)\right).
$$
This can be shown to be zero unless $j = \ell(w)$, in which case it is the trivial module $k$.   One then needs to know the vanishing of certain higher right derived induction functors:
\vskip.2cm
\noindent
{\em Assumption 1:} $R^i\ind_{P_J}^G S^{\bullet}(\ul_J^*) = 0$ for $i > 0$.
\vskip.2cm
\noindent
If true, one may conclude that
\begin{itemize}
\item $\opH^{2\bullet + 1}(G_1,k) = 0$
\item $\opH^{2\bullet}(G_1,k)^{(-1)} \cong \ind_{P_j}^GS^{\bullet}(\ul_J^*)$.
\end{itemize}

One would like to identify $\ind_{P_j}^GS^{\bullet}(\ul_J^*)$ with  $k[G\cdot \ul_J]$. That can be shown to hold under a second geometric assumption:
\vskip.2cm
\noindent
{\em Assumption 2:} The Richardson orbit closure $G\cdot \ul_J$ is a normal subvariety of ${\mathcal N}$.
\vskip.2cm
\noindent

Note that Assumptions 1 and 2 trivially hold in the case $p = h$, but that occurs only in type $A_n$ where $p = n + 1$, which is precisely one of our excluded cases.   This will be discussed more below.   The assumptions are also known to hold in the case $p = h - 1$, where $G\cdot \ul_J$ is the closure of the subregular orbit.  More generally, with the assumptions on $p$ as above, cohomology will arise only in even degrees and one may identify $\opH^{\bullet}(G_1,k)^{(-1)}$ with $k[G\cdot \ul_J]$ under the following assumptions:

\vskip.2cm
\begin{center}
\begin{tabular}{|c|c|c|c|c|c|c|}
\hline
Type & $A_n$, $C_n$, $D_n$ & $B_n$ & $E_6$, $G_2$ & $F_4$ & $E_7$ & $E_8$\\
\hline
Condition & $p > h/2$ & $p \geq h/2$ & $p \geq 5$ & $p \geq 7$ & $p \geq 11$ & $p \geq 17$\\
\hline
\end{tabular}.
\end{center}

\vskip.2cm
Again, the reader is referred to \cite{BNPP} for more details and explicit information on the subset $J$ (and element $w \in W$).   There one will also find a discussion of the type $A_n$ case when $p$ divides $n + 1$, where the answer has a conjecturally different format (as already seen in \cite{AJ}).

\subsection{$B_1$-cohomology}\label{S:B1} As part of the computation of $\opH^{\bullet}(G_1,k)$, when $p > h$, one finds that $\opH^{\bullet}(B_1,k) \cong S^{\bullet}(\ul^*)^{(1)}$ (with the generators in degree 2); a result that fails when $p \leq h$.   By the nature of the small prime argument just discussed, no new  information is gained on the $B_1$-cohomology and this remains an open question in general.  There are some known results for special linear groups as discussed in the next subsection and more generally in low degrees as discussed in Section \ref{S:smalldeg}.

\subsection{Higher Frobenius kernels}\label{S:grk} Much less is known about $\opH^{\bullet}(G_r,k)$ for $r > 1$.   For $G = SL_2$, $\opH^{\bullet}(U_r,k)$ was identified in \cite{AJ}, from which some observations on $B_r$-cohomology were made for $r = 1, 2$.   Ngo \cite{N13} expanded on this, providing complete calculations of the $B_r$- and $G_r$-cohomology.  This is done in the more general context (discussed further in Section \ref{S:nontriv}) of computations of $\opH^i(B_r,\la)$ and $\opH^i(G_r,H^0(\la))$ for $\la \in X^+$.    In this case, one has an affirmative answer to Question \ref{Q:ind}.   For the rings $\opH^{\bullet}(B_r,k)$ and $\opH^{\bullet}(G_r,k)$, Ngo also identified the ring structure of the reduced rings.  

More recent work of Ngo \cite{N19} for arbitrary $G$ shows that $B_r$-cohomology is related to $B_1$-cohomology.   More precisely, it is shown that, 
for $i < \frac{p}{c}$, 
$$
\opH^i(B_r,k) \cong \opH^i(B_1,k)^{(r-1)} \cong \begin{cases}
S^{\frac{i}{2}}(\ul^*)^{(r)} &\text{ if } i \text{ is even},\\
0 &\text{ if } i \text{ is odd},
\end{cases}
$$
where $c$ depends on the root system and is given by the following table:
\vskip.2cm
\begin{center}
\begin{tabular}{|c|c|c|c|c|c|c|c|}
\hline
Type & $A_n$ & $B_n$, $C_n$, $D_n$ & $E_6$ & $E_7$ & $E_8$ & $F_4$ & $G_2$\\
\hline
$c$ & 1 & 2 & 3 & 4 & 6 & 4 & 3\\
\hline
\end{tabular}.
\end{center}
\vskip.2cm\noindent
Using the ideas of Section \ref{S:ind} one may conclude that $\opH^i(G_r,k)^{(-r)} \cong \ind_B^G\left(\opH^i(B_r,k)^{(-r)}\right)$ for $i < \frac{p}{c}$.

Such a result had been foreshadowed by earlier observations in type $A_n$ made by Friedlander and Parshall \cite{FP862} and Kaneda, Shimada, Tezuka, and Yagita \cite{KSTY} where such a relationship was shown to hold in degrees less than $2p - 1$ for $p > h$.   Kaneda et al.~explicitly computed $\opH^{2p-1}(B_r,k)$ for $r > 1$, finding it to be non-zero, thus demonstrating that the $2p-1$ bound was strict.   Note also that $\opH^{2p-1}(G_r,k)^{(-r)} \cong H^0(\al_0)$.  The work of Kaneda et al.~also gave a complete description of $\opH^{\bullet}(B_2,k)$ for $SL_3$ for $p \geq 3$.

More recently, Friedlander \cite{F19} has investigated $U_r$-cohomology, specifically considering the case that $U$ is the unipotent radical of $GL_3$. The cohomology ring $\opH^{\bullet}(U_r,k)$ is approximately identified in a sense.  It is shown that there is an algebra given by known generators and relations that embeds into $\opH^{\bullet}(U_r,k)$ with image containing all $p$th powers.  This arises as a special case of a more general theory developed for the unipotent radical $U_J$ of a standard parabolic in the general linear group $GL_n$ and considering the quotient by a term in its descending central series.


\section{Cohomology of Frobenius Kernels with Non-trivial Coefficients}\label{S:nontriv}

\subsection{} More generally, one may consider cohomology groups $\opH^i(G_r,M)$, $\opH^i(B_r,M)$, or $\opH^i(U_r,M)$ for a $G_r$-module (or $B_r$-module or $U_r$-module) $M$. Over $G_r$, natural modules of interest are a simple module $L(\la)$ or a standard induced module $H^0(\la)$.   Thanks to relationships (seen above) between $B_r$- and $G_r$-cohomology, much more is known in the latter case.

\subsection{Induced modules}\label{S:AJ} The seminal result in this context (for $r = 1$) is due to Andersen and Jantzen \cite{AJ}, with certain exceptional cases dealt with later by Kumar, Lauritzen, and Thomsen \cite{KLT}.  Assume $p > h$ and $\la \in X^+$. Then
\begin{align*}
\opH^i(G_1,H^0(\la))^{(-1)} &\cong \ind_{B}^{G}\left(\opH^i(B_1,\la)^{(-1)}\right)\\
&\cong
\begin{cases}
\ind_B^G\left(S^{\frac{i - \ell(w)}{2}}(\ul^*)\otimes\mu\right) &\text{ if } \la = w\cdot 0 + p\mu,\\
0 & \text{ otherwise,}
\end{cases}
\end{align*}
where $w \in W$ and $\mu \in X^+$.  Implicit in the statement is not only a computation of $\opH^i(G_1,H^0(\la))$ but also the computation of $\opH^i(B_1,\la)$ generalizing that of Section \ref{S:B1} in the case $\la = 0$.   Again, as noted in \cite{AJ}, this formula breaks down for $p \leq h$. 

For small primes and higher $r$, there is only one general result, and that is when $G = SL_2$, found in the aforementioned work of Ngo \cite{N13}.    Indeed, the results mentioned in Section \ref{S:grk} are special cases of this work, where explicit computations of $\opH^i(U_r,\la)$, $\opH^i(B_r,\la)$, and $\opH^i(G_r,H^0(\la))$ are given for all $i$, $r$, and $\la \in X^+$.

\subsection{Small degrees}\label{S:smalldeg} In degree 1, Jantzen \cite{Jan91} computed all $\opH^1(B_1,\la)$ and $\opH^1(G_1,H^0(\la))$.  Those computations were extended inductively to $B_r$ and $G_r$ by Bendel, Nakano, and Pillen \cite{BNP04} and then to degree 2  for $p \geq 3$ \cite{BNP07}.  The degree 2 and $p = 2$ case was done by Wright \cite{W}.   Bendel, Nakano, and Pillen \cite{BNP16} extended such computations to degree 3 as long as the prime is not too small: $p \geq 5$ in Types $A_n$ ($p \neq 5$ in type $A_4$ and $p\neq 7$ in type $A_6$), $C_n$, $D_n$ or $E_n$; and  $p \geq 7$ in types $B_n$ ($n \geq 3$), $F_4$, or $G_2$. 

The general approach follows the ideas of Section \ref{S:trivial} for trivial coefficients, but becomes more intricate as the degree increases.   Following Section \ref{S:ind} the goal is to conclude that
$$
\opH^i(G_r,H^0(\la))^{(-1)} \cong \ind_{B}^{G}\left(\opH^i(B_r,\la)^{(-1)}\right),
$$
which requires not simply knowledge of $\opH^i(B_r,\la)^{(-1)}$  per se, but the specific nature of those groups to be able to conclude that higher right derived induction functors vanish.    $B_r$-cohomology is related to $B_1$-cohomology through the Lyndon-Hochschild-Serre spectral sequence.  Then $B_1$-cohomology is related to $U_1$-cohomology:
$$\opH^i(B_1,\la) \cong \left(\opH^i(U_1,\la)^{T_1}\right) \cong \left(\opH^i(U_1,k)\otimes\la\right)^{T_1}.$$
Lastly, one  uses \eqref{E:AJss} to relate $U_1$-cohomology to $\ul$-cohomology.\footnote{For $p = 2$, the spectral sequence is not available.  There is a complex that may be used to compute $U_1$-cohomology (cf. \cite[Lemma I.9.15]{Jan}).}   

In practice, the computational process is more intertwined than the above overview might suggest. For example, to make computations of $\ul$-cohomology one sometimes has to partially reverse the above process: via the connections between $\ul$, $U_1$, and $B_1$, using the Lyndon-Hochschild-Serre spectral sequence for $B_1$ in $B$, one may use information about $B$-cohomology to get $\ul$-cohomology information.  In degree 3, one also begins to see that not only is information on $\opH^i(\ul,k)$ necessary, but also Lie algebra cohomology for non-trivial coefficients, for example, $\opH^i(\ul,\ul^*)$.   In higher degrees, this approach would potentially need $\opH^i(\ul,S^j(\ul^*))$.  

An interesting feature of this spiraling/inductive strategy is that while the initial aim might be understanding $\opH^i(G_r,H^0(\la))$, one obtains new cohomology results for $B_r$, $B$, and $\ul$ that are of interest in their own right.   The latter two cases will be discussed in Sections \ref{S:alg} and \ref{S:Lie} respectively.  

If the underlying root system $\Phi$ is of type $C_2$ or $F_4$ with $p = 2$ or $G_2$ with $p = 3$, $G$ admits a purely inseparable isogeny whose square is the Frobenius morphism.  This gives rise to so-called half Frobenius kernels $G_{r/2}$ for an odd positive integer $r$.\footnote{These isogenies are also involved in the construction of the finite Suzuki-Ree groups.}  Similar to the $B_r$- and $G_r$-computations above, Radu \cite{Rad} has computed $\opH^1(B_{r/2},\la)$ and $\opH^1(G_{r/2},H^0(\la))$ in these cases.

\subsection{Lie algebra cohomology}\label{S:Lie}  For a Lie algebra $\mathfrak{h}$ the ordinary Lie algebra cohomology is equivalently that of its universal enveloping algebra. For an $\mathfrak{h}$-module $M$, $\opH^i(\mathfrak{h},M)$ may be computed from a complex of exterior powers: $M\otimes\Lambda^{\bullet}(\mathfrak{h}^*)$.  See for example \cite[I.9.17]{Jan} and note that $\opH^i(\mathfrak{h},M) = 0$ for $i > \dim\mathfrak{h}$.  The primary case of interest here is for the Lie algebra $\ul$ with trivial coefficients.

A classical result of Kostant \cite{Kos}\footnote{Kostant's result is actually for more general coefficients in a simple module, as are subsequent prime characteristic results that are mentioned below.} is
$$
\opH^i(\ul,k) \cong \bigoplus_{\ell(w) = i, w\in W} -w\cdot 0.
$$ 
This was extended to prime characteristic by Friedlander and Parshall \cite{FP862} for $p > h$ and improved to $p \geq h - 1$ by Polo and Tilouine \cite{PT} and, more recently, by the University of Georgia Vigre Algebra group \cite{UGA}.  There are also more general versions of this result for $\ul_J$ associated to a parabolic subgroup for a set of simple roots $J$.  

For $p < h - 1$, one will necessarily have cohomology classes with weights as above, but there may be more.  Indeed, the University of Georgia VIGRE Algebra Group showed that there is always more cohomology in a global character sense:
$$
\operatorname{ch}\opH^{\bullet}(\ul,k) \neq \operatorname{ch}\opH^{\bullet}(\ul,{\mathbb C}).
$$
As a result, there is necessarily some degree $i$ where $\dim\opH^i(\ul,k) > \dim\opH^i(\ul, \mathbb{C})$, but that does not mean that the dimensions differ in all degrees.  As part of the University of Georgia VIGRE Algebra Group work, there is a \textsc{Magma} program written by Brian Boe to compute $\opH^i(\ul,k)$.  The author also had undergraduate students at the University of Wisconsin-Stout making computer computations.  The dimensions of all $\opH^i(\ul,k)$ are known at least for all rank 2 and 3 groups, as well as in type $A_4$ and $A_5$.   Partial results are known in some higher ranks, with practical computational power being the limiting factor.

For small primes, as implied in Section \ref{S:smalldeg}, complete degree 1 computations were made by Jantzen \cite{Jan91} and degree 2 computations by Bendel, Nakano, and Pillen \cite{BNP04} and Wright \cite{W}.  One finds that Kostant's Theorem holds in all types if $p > 3$.  In degree 3, Bendel, Nakano, and Pillen \cite{BNP16} showed that Kostant's Theorem holds with stronger conditions on $p$ as given in the following chart:

\vskip.2cm
\begin{center}
\begin{tabular}{|c|c|c|c|}
\hline
Root System & $A_n$, $n \leq 3$; $B_2$ & $A_n$, $n \geq 4$; $B_3$; $C_n$, $n \geq 3$; $D_n$; $E_n$; $G_2$ & $B_n$, $n \geq 4$; $F_4$\\
\hline
$p \geq $ & 3 & 5 & 7\\
\hline
\end{tabular}.
\end{center}

\subsection{Simple modules} Returning to the question of $G_1$- or $G_r$-cohomology, a more challenging problem, and one where much less is known, is to compute $\opH^i(G_1,L(\la))$ for $\la \in X_1$ when $L(\la) \neq H^0(\la)$. One may try to use the long exact sequence in cohomology associated to the short exact sequence
$$
0 \to L(\la) \to H^0(\la) \to H^0(\la)/L(\la) \to 0.
$$
Some results in that direction were provided in degree 1 by Jantzen \cite{Jan91} for $G_1$ and by Bendel, Nakano, and Pillen \cite{BNP04} for $r > 1$.   This is a special case of determining $\Ext^1_{G_1}(L(\la),L(\mu))$ for $\la,\mu \in X_1$, where some results are known in small rank: for example, Yehia \cite{Yeh} for type $A_2$, Sin \cite{Sin} and Dowd and Sin \cite{DS} in characteristic 2 for ranks at most 4, Sin \cite{Sin} for type $G_2$ and $p = 3$, and Lin \cite{Lin} for type $G_2$ and $p \geq 3(h-1)$.  Radu \cite{Rad} has recently made computations of $\Ext^1_{G_r}(L(\la),L(\mu))$ for $\la,\mu \in X_r$ for general $r$ in types $C_2$ or $F_4$ with $p = 2$ and $G_2$ with $p = 3$, as well as for the half Frobenius kernel $G_{r/2}$ mentioned in Section \ref{S:smalldeg}.

A finite-dimensional rational $G$-module $M$ is said to be a {\em tilting} module if it admits both a good and Weyl filtration (cf. \cite[II.E.1, II.4.16, II.4.19]{Jan}).  That is, a filtration with factors of the form $H^0(\la)$ (good filtration) or $V(\la)$ (Weyl filtration) for $\la \in X^+$.  A fundamental question for Ext-groups is the following.

\begin{quest}\label{Q:ExtTilt} Let $\la, \mu \in X_1$. Is $\Ext_{G_1}^1(L(\la),L(\mu))^{(-1)}$ a tilting module?
\end{quest}

This is known to hold for large primes, where it follows from a stronger statement.    For $p \geq 3h - 3$, Andersen \cite[Prop. 5.5]{And} showed that $\Ext_{G_1}^1(L(\la),L(\mu))^{(-1)}$ is semi-simple with each composition factor $L(\si)$ having highest weight $\si$ lying in the closure of the fundamental alcove.  For such a $\si$, one has $L(\si) = H^0(\si) = V(\si)$, from which the existence of good and Weyl filtrations are immediate.   The bound on the prime for this stronger statement was lowered to $p \geq 2h - 2$ by Bendel, Nakano, and Pillen \cite[Cor. 5.5B]{BNP04} and then to $p \geq 2h - 4$ by Bendel, Nakano, Pillen, and Sobaje \cite[Thm. 4.3.1]{BNPS23}.  

For small primes, there exist negative answers to Question \ref{Q:ExtTilt}.   For example, when $p = 2$ and $\Phi$ is of type $B_n$ ($n \geq 3$) or $G_2$, Jantzen \cite[Prop. 6.9]{Jan91} showed that $\opH^1(G_1,L(\varpi_2))^{(-1)} \cong \Ext_{G_1}^1(k,L(\varpi_2))^{(-1)} \cong H^0(\varpi_1)$ (following the Bourbaki ordering of roots).  But $H^0(\varpi_1)$ is not simple (and hence not tilting).   For $p = 3$ and $\Phi$ being of type $C_3$, from \cite[Prop. 4.1, \S 4.2, Prop. 4.5]{Jan91}, one can show that $\Ext_{G_1}^{1}(k,L(\varpi_1 + \varpi_2 + \varpi_3))^{(-1)}$ is not tilting (cf. \cite[\S 4.5.3]{BNPS20}).   Recent work of Bendel, Nakano, Pillen, and Sobaje \cite{BNPS20, BNPS22, BNPS23} has shown that there is close connection between Question \ref{Q:ExtTilt} and a conjecture of Donkin on tilting modules.   Let $\la \in X_r$ and $Q_r(\la)$ denote the injective hull (equivalently, projective cover) of $L(\la)$ as a $G_r$-module.   A yet unresolved question dating to work of Humphreys and Verma \cite{HV} is whether $Q_r(\la)$ admits the structure of a rational $G$-module.\footnote{It is known to hold for $p \geq 2h - 4$, as well as for smaller primes for some small rank groups.}  Donkin \cite{Don} conjectured that this lift should be a specific indecomposable tilting module: $T(2(p-1)\rho + w_0\la)$.   The negative answers to Question \ref{Q:ExtTilt} noted above have led to the first known counterexamples to Donkin's conjecture (cf. \cite{BNPS20, BNPS22}).


\section{$B$-cohomology}\label{S:alg}

\subsection{}\label{S:Bintro} Leaving the domain of Frobenius kernels and returning to the full algebraic group setting, consider the Borel subgroup $B$.   Its cohomology with trivial coefficients agrees with that of $G$ and is well-known:
$$
\opH^i(G,k) \cong \opH^i(B,k) = \begin{cases}
k &\text{ if } i = 0,\\
0 &\text{ if } i > 0.
\end{cases}
$$
Indeed, for a rational $G$-module $M$, one has (cf. \cite[Cor. II.4.7]{Jan}) $\opH^i(G,M) \cong \opH^i(B,M)$ and one might attempt to use $B$-information to translate to $G$.  

However, one is often interested in $\opH^i(B,M)$ where $M$ is a $B$-module that may not be a $G$-module.  In general $B$-cohomology can be identified with an inverse limit over $r$ of $B_r$-cohomology (cf. \cite[Cor. II.4.12]{Jan}):
\begin{equation}\label{E:BrB}
\opH^i(B,M) \cong \underset{\longleftarrow}{\lim}\opH^i(B_r,M);
\end{equation}
a statement that also hods for $G$ or any parabolic subgroup thereof.  
One also has some general knowledge on the cohomological degrees in which $\opH^i(B,M)$ can be non-zero based on weights of the module $M$ (cf. \cite[Prop. II.4.10]{Jan}), but specific computations are limited.  One fundamental case of interest (suggested by earlier discussions) is a one-dimensional $B$-module.

\subsection{One-dimensional $B$-modules}\label{S:Bla} In characteristic zero, it is well-known that, for $\la \in X$, $\opH^i(B,\la)$ is zero unless $\la = w\cdot 0$. For such a $\la$, non-vanishing occurs only in degree $i = \ell(w)$, where one has $\opH^{\ell(w)}(B,\la) = k$.  This follows from the Borel-Bott-Weil Theorem that describes $R^i\ind_B^G\la$ for all $\la$ and $i$, along with the fact that $G$-cohomology vanishes in positive degree. 

In prime characteristic, much less is known about $\opH^i(B,\la)$.   In degree 1, Andersen \cite{And} showed that the only non-zero cohomology is the following:
$$
\opH^1(B,-p^m\al) = k
$$
for $m \geq 0$, $\al \in S$.

The answer is also known in degree 2 due to work of Andersen \cite{And}, O'Halloran \cite{OHal}, Bendel, Nakano, and Pillen \cite{BNP07}, and Wright \cite{W}. The cohomology is always at most one-dimensional.  For $p > 3$, $\opH^2(B,\la)$ is non-zero in the following cases:
$$
\la = \begin{cases}
p^iw\cdot 0, \quad &\text{ for }i \geq 0, \ell(w) = 2,\\
-p^i\al, \quad &\text{ for } i > 0, \al \in S,\\
-p^i\al - p^j\be, &\text{ for } 0 \leq i < j, \al, \be \in S.
\end{cases}
$$
For $p = 3$, additional cases occur in type $G_2$, and, for $p = 2$, additional cases occur in all non-simply laced root systems (i.e., types $B_n$, $C_n$, $F_4$, and $G_2$).

There are many more cases where $\opH^3(B,\la) \neq 0$, as well as cases where the cohomology is 2-dimensional (even for large $p$).\footnote{For example, one such case is if $\la = -p^i\al - p^j\be$ for $i > j \geq 1$ and $\al, \be \in S$.}  For $p > h$, complete computations of $\opH^3(B,\la)$ were given by Andersen and Rian \cite{AR} using the Lyndon-Hochschild-Serre spectral sequence for $B_1$ in $B$, along with known computations for $B_1$-cohomology (as discussed in Section \ref{S:smalldeg}).    Bendel, Nakano, and Pillen \cite{BNP16} recovered the \cite{AR} results and extended them to smaller primes (with restrictions on $p$ as given in the table in Section \ref{S:Lie} by using \eqref{E:BrB} and computations of $\opH^i(B_r,\la)$.

\subsection{$U$-cohomology} While the computation of $\opH^i(B,k)$ is trivial, the closely related question of determining $\opH^i(U,k)$ is more complex and seems to have received minimal study.  Fairly recently, Friedlander \cite{F19} initiated an attempt to better understand $\opH^i(U,k)$.  In particular, he observed that $\opH^i(U,k)$ can be identified with the inverse limit over $r$ of the $\opH^i(U_r,k)$, analogous to the long-known result \eqref{E:BrB}.  If computations could be made for $U_r$, analogous to the approached just mentioned for $B$-cohomology, this would provide an avenue for computing $U$-cohomology.


\section{Rational and Generic Cohomology}\label{S:rat}

\subsection{} One of the most fundamental connections between cohomology groups is the seminal work of Cline, Parshall, Scott, and van der Kallen \cite{CPSvdK} relating $G$-cohomology and $\gfq$-cohomology. Specifically, they showed that, for a finite-dimensional rational $G$-module $M$ and sufficiently large $r$ and $s$ (depending on $i$), the restriction map
$$
\opH^i(G,M^{(s)}) \to \opH^i(\gfq,M^{(s)})
$$
is an isomorphism.   Since the Frobenius morphism on $G$ is an automorphism on the subgroup $\gfq$ (for any $r$), $M^{(s)} \cong M$ as a $\gfq$-module for any $s \geq 1$. In other words, the above isomorphism becomes
\begin{equation}\label{E:gfqtoG}
\opH^i(G,M^{(s)}) \cong \opH^i(\gfq,M),
\end{equation}
noting that the right-hand side is independent of $s$.  

From the perspective of $\gfq$-cohomology, \eqref{E:gfqtoG} says that one can equate $\opH^i(\gfq,M)$ with $G$-cohomology as long as one applies a sufficiently high twist to $M$.  Recent developments in regards to removing (in a sense) the need for the twist are discussed in Section \ref{S:shift} below.   The isomorphism \eqref{E:gfqtoG} leads to two stability consequences:

\vskip.2cm\noindent
{\em Rational Stability:} The cohomology $\opH^i(G,M^{(s)})$ eventually stabilizes as $s$ increases. That is, for a given $i$, there exists an $s$ such that 
\begin{equation}\label{E:rat}
\opH^i(G,M^{(s)}) \cong \opH^i(G,M^{(s+1)}).
\end{equation}

\vskip.2cm\noindent
{\em Generic Cohomology:} The cohomology $\opH^i(\gfq,M)$ eventually stabilizes as $r$ increases.   This stable value is known as the {\em generic} cohomology $\opH^i_{\rm{gen}}(G,M)$.

\vskip.2cm
A new perspective on these ideas was taken in work of Bendel, Nakano, and Pillen \cite{BNP14}.  This will be discussed in Sections \ref{S:genrev} and \ref{S:stab} which depends heavily on certain induced modules that are considered in the next subsection.

\subsection{Inducing from $\gfq$ to $G$}\label{S:indgfq} A direct relationship between $\gfq$-cohomology and $G$-cohomology is given by generalized Frobenius reciprocity (or Shapiro's lemma). For a $\gfq$-module $M$, one has
\begin{equation}\label{E:ind}
\opH^i(\gfq,M) \simeq \opH^i(G,\ind_{\gfq}^G(M)).
\end{equation}
For a rational $G$-module $M$,  by the tensor identity, one has $\ind_{\gfq}^G(M) \cong M\otimes \ind_{\gfq}^G(k)$ and hence
$$\opH^i(\gfq,M) \simeq \opH^i(G,M\otimes\ind_{\gfq}^G(k)).$$
Unfortunately, the module $\ind_{\gfq}^G(k)$ is infinite-dimensional.

As observed by Bendel, Nakano, and Pillen \cite[Prop. 2.4]{BNP11} (cf.~also \cite{PSS, BNPPSS}), $\ind_{\gfq}^G(k)$ admits a filtration with factors of the form $H^0(\lambda)\otimes H^0(\lambda^*)^{(r)}$, where there is precisely one such factor for each $\la \in X^+$.  That is, there exists a filtration
$$
0 \subseteq F_0 \subseteq F_1 \subseteq F_2 \subseteq \cdots \subseteq \ind_{\gfq}^G(k),
$$
where $F_i/F_{i - 1} \cong H^0(\la_i)\otimes H^0(\la_i^*)^{(r)}$ for $\la_i \in X^+$, $\cup_i F_i = \ind_{\gfq}^G(k)$, and each $\la \in X^+$ arises precisely once as a $\la_i$. Moreover, this filtration behaves nicely with respect to the size of weights.   Given a positive integer $m$, $\ind_{\gfq}^G(k)$ admits a finite-dimensional submodule $I_m$ that has a (necessarily finite) filtration as above, where each $\la_i$ appearing in a filtration factor satisfies $\langle \la_i, \al_0^{\vee}\rangle \leq m$.

Such a filtration had been foreshadowed in earlier work of Bendel, Nakano, and Pillen \cite{BNP01, BNP02, BNP040, BNP041, BNP06}, where truncations of $\ind_{\gfq}^G(k)$ (to certain bounded weight categories) were considered and used to directly relate $\gfq$- and $G$-cohomology (and extensions more generally).    The reader is referred to \cite[\S 5]{Nak} for an overview of those results.

\subsection{Generic cohomology revisited}\label{S:genrev}  The filtration of Section \ref{S:indgfq} was used in \cite{BNP14} to give a new proof of \eqref{E:gfqtoG}.   The idea is to to consider a short exact sequence
$$
0 \to k \to \ind_{\gfq}^G(k) \to N \to 0.
$$
Tensoring with $M^{(s)}$ gives another exact sequence:
\begin{equation}\label{E:exact}
0 \to M^{(s)} \to M^{(s)}\otimes\ind_{\gfq}^G(k) \to M^{(s)}\otimes N \to 0.
\end{equation}
One can then take the long exact sequence in $G$-cohomology associated to \eqref{E:exact}. Using the identification $\opH^i(G,M^{(s)}\otimes \ind_{\gfq}^G(k)) \cong \opH^i(G,\ind_{\gfq}^G(M^{(s)})) \cong \opH^i(\gfq,M^{(s)})$, this becomes
$$
\begin{array}{cclclclcc}
0 & \longrightarrow & \operatorname{Hom}_G(k,M^{(s)}) & \stackrel{\res}{\longrightarrow} & \operatorname{Hom}_{G({\mathbb F}_{q})}(k,M^{(s)}) &\longrightarrow & \operatorname{Hom}_G(k,M^{(s)}\otimes N) &&\\
& \longrightarrow & \operatorname{H}^1(G,M^{(s)}) & \stackrel{\res}{\longrightarrow} & \operatorname{H}^1(G({\mathbb F}_{q}),M^{(s)}) & \longrightarrow & \operatorname{H}^{1}(G,M^{(s)}\otimes N) &&\\
& \longrightarrow & \operatorname{H}^2(G,M^{(s)}) & \stackrel{\res}{\longrightarrow} & \operatorname{H}^2(G({\mathbb F}_{q}),M^{(s)}) & \longrightarrow & \operatorname{H}^{2}(G,M^{(s)}\otimes N) &&\\
& \longrightarrow & \operatorname{H}^3(G,M^{(s)}) & \stackrel{\res}{\longrightarrow} & \operatorname{H}^3(G({\mathbb F}_{q}),M^{(s)}) & \longrightarrow & \operatorname{H}^{3}(G,M^{(s)}\otimes N) &\longrightarrow &\cdots.  
\end{array}
$$
One sees that the groups in the third column are the obstruction to the restriction map being an isomorphism.  To show that $\opH^{i}(G,M^{(s)}\otimes N)$ vanishes, one shows the vanishing of all $\opH^i(G,M^{(s)}\otimes H^0(\la_j)\otimes H^0(\la_j^*)^{(r)})$ for the $H^0(\la_j)\otimes H^0(\la_j^*)^{(r)}$ appearing as filtration factors in $N$.

While these generic cohomology statements are existence in nature, both the original work \cite{CPSvdK} and the newer work \cite{BNP14} give explicit values for $r$ and $s$ where \eqref{E:gfqtoG} is guaranteed to hold.    One distinction between the methods is that the values in \cite{CPSvdK} depend on the root system (in addition to the degree $i$ and prime $p$), whereas those in \cite{BNP14} depend only on the degree and the prime.   In \cite[\S 8]{BNP14} a comparison of bounds is made, showing a general improvement (decrease) in the size of $r$ and $s$ needed.

\subsection{Rational stability revisited}\label{S:stab} As noted above, a consequence of \eqref{E:gfqtoG} is the rational stability result \ref{E:rat}.  In \cite{BNP14}, the authors gave an independent proof of this using new results on $B$-cohomology of one-dimensional $B$-modules as discussed in Section \ref{S:Bla}.

\subsection{Shifted cohomology}\label{S:shift}  A particular case of interest is when the rational $G$-module $M$ is a simple module $L(\la)$ with $\la \in X_r$, so that $L(\la)$ remains simple upon restriction to $\gfq$. Then one has $\opH^i(\gfq,L(\la)) \cong \opH^i(G,L(\la)^{(s)})$.  It turns out that one can remove the $s$-twist in a qualified sense.  The potential for such a result was seen for example in \cite[Thm. 5.5]{BNP06}, where it was shown that, for sufficiently large $r$ (depending on the prime $p$ and Coxeter number $h$), either $\opH^1(\gfq,L(\la)) \cong \opH^1(G,L(\la))$ or $\opH^1(\gfq,L(\la)) \cong \opH^1(G,L(\tilde{\la}))$ for a {\em different} $p^r$-restricted weight $\tilde{\lambda}$ (obtained by a certain shifting algorithm applied to the weight $\la$).   A version for extensions between two simple modules was also shown to hold, where the highest weights of either or both simple modules might need to be shifted.  

Later work of Stewart \cite{Ste2} for degree 2 cohomology with $G$ being $SL_3$ suggested that such a result might hold in higher degrees.  Parshall, Scott, and Stewart \cite{PSS} showed that one may choose $r$ sufficiently large (depending on the root system and cohomological degree $m$) so that, for $i \leq m$ and $\la \in X_r$, $\opH^i(\gfq,L(\la)) \cong \opH^i(G,L(\la'))$ for some potentially shifted weight $\la' \in X_r$.   A version for Ext-groups between simples was also given.


\section{Vanishing Ranges}\label{S:van}

\subsection{} Dating back to work of Quillen \cite{Q}, it has been observed that $\opH^i(\gfq,k)$ is zero in ``low'' positive degrees.  That is
$$
\opH^i(\gfq,k) = 0 \text{ for } 0 < i < m
$$
for some $m$ depending on $G$ and $r$.  Ideas discussed in the previous section have been used to provide new results on this problem, discussion of which may also be found in \cite[\S 6]{Nak}. 

\subsection{The general linear group} Quillen \cite{Q} showed that $\opH^i(GL_n({\mathbb F}_{q}),k) = 0$  for $0 < i < r(p - 1)$, while making no claim about what happened in degree $r(p-1)$.  The question of interest then is to identify precisely the least positive degree with non-zero cohomology.  Friedlander and Parshall \cite{FP83} (for $p > 2$) expanded Quillen's original vanishing range to $0 < i < r(2p - 3)$ by showing such vanishing for the Borel subgroup $B(\fq)$. Further, they showed the bound was sharp for $B(\fq)$. That is, $\opH^{r(2p - 3)}(B({\mathbb F}_{q}),k) \neq 0$.  But they were not able to show non-vanishing of $\opH^{r(2p-3)}(GL_n(\fq),k)$.   Barbu \cite{B} showed that $\opH^{r(2p - 2)}(GL_n({\mathbb F}_{q}),k) \neq 0$ and conjectured that it should be related (via a Bockstein map) to a non-zero class in degree $r(2p-3)$.   Non-vanishing in degree $r(2p-3)$ was shown by Bendel, Nakano, and Pillen \cite{BNP11} under the condition that $p \geq n + 2$. The bound on the prime was lowered to $p \geq n$ by Sprehn \cite{S} who took a more constructive approach.  More recent work of Lahtinen and Sprehn \cite{LS} constructed non-zero cohomology elements when $p < n$.  See more discussion in Section \ref{S:smallvan} below.

\subsection{Other groups} Quillen \cite{Q} commented (without any specifics) that there should be a similar vanishing range for other classical groups.   Vanishing ranges (but not the least degree of non-vanishing) were found by Friedlander \cite{F76} for special orthogonal and symplectic groups and by Hiller \cite{H} for simply connected $G$ in all types.   

The first precise bounds were given in work of Bendel, Nakano, and Pillen \cite{BNP11, BNP12}.   As with earlier arguments for Frobenius kernels, one may focus on the $r = 1$ case (i.e., $\gfp$) and then attempt to use induction to obtain results for higher $r$.   One has the isomorphism
$$\opH^i(\gfp,k) \simeq \opH^i(G,\ind_{\gfp}^G(k))$$ 
of \eqref{E:ind} along with the filtration of $\ind_{\gfp}^G(k)$ discussed in Section \ref{S:indgfq}.  One then wishes to show the vanishing of 
$$\opH^i(G,H^0(\lambda)\otimes H^0(\lambda^*)^{(1)}) \cong \Ext^i_G(V(\la)^{(1)},H^0(\la))$$
in small degrees.  To that end, one may use the Lyndon-Hochschild-Serre spectral sequence for $G_1 \unlhd G$:
$$
E_2^{i,j} = \Ext^i_{G/G_1}(V(\la)^{(1)},\opH^j(G_1,H^0(\la))) \Rightarrow \Ext^{i+j}_{G}(V(\la)^{(1)},H^0(\la)).
$$
One may now apply the result of Andersen-Jantzen from Section \ref{S:AJ} on $\opH^j(G_1,H^0(\la))$ (where the assumption $p > h$ is needed).   This leads to combinatorial questions based on weight spaces in $S^{\bullet}(\ul^*)$ and Kostant's partition function.  

For $G$ simple, simply connected as in our standing assumption, for $p > h$, precise bounds are obtained for $\gfp$ in all types except $F_4$ and $G_2$.    Generically, $\opH^i(\gfp,k) = 0$ for $0 < i < 2p-3$ and $\opH^{2p-3}(\gfp,k) \neq 0$.  However, that is not universally the case with exceptions in type $C_n$, ``small'' rank, and in type $A_n$ when $p = n + 2 = h + 1$.   A complete summary of these bounds is given in the following table, where the sharp bound $D$ means that $\opH^D(\gfp,k) \neq 0$ and $\opH^i(\gfp,k) = 0$ for $0 < i < D$.

\vskip.2cm
\begin{center}
\begin{tabular}{|c|c|c|}
\hline
Root System &  $p$ & Sharp Vanishing Bound $D$ \\
\hline
$A_n, n \geq 4$ & $p > n + 2$  & $2p - 3$\\
\hline
$A_n, n \geq 3$ & $p = n + 2$ & $p - 2$ \\
\hline
$A_3$ & $p > 5$ & $2p - 6$  \\
\hline
$A_2$ & $p = 3k + 1$, $k \geq 2$ & $2p - 6$ \\
\hline
$A_2$ & $p = 3k + 2$, $k \geq 1$ & $2p - 3$ \\
\hline
$B_n, n \geq 7$ & $p > 2n$ & $2p - 3$  \\
\hline
$B_n, n \in \{5, 6\}$ & $p > 13$ & $2p - 3$  \\
\hline
$B_n, n \in \{5,6\}$ & $p = 13$ & $2p - 5$  \\
\hline
$B_5$ & $p = 11$ & $2p - 7 = 3$  \\
\hline
$B_n$, $n \in \{3,4\}$ & $p > 2n$ & $2p - 8$  \\
\hline
$C_n, n \geq 1$ & $p > 2n$ & $p - 2$  \\
\hline
$D_n, n \geq 4$ & $p > 2n - 2$ & $2p - 2n$ \\
\hline
$E_6$ & $p > 13$ & $2p - 3$ \\
\hline
$E_6$ & $p = 13$ & $2p - 10 = 16$ \\
\hline
$E_7$ & $p > 23$ & $2p - 3$ \\
\hline 
$E_7$ & $p = 23$ & $2p - 7 = 39$ \\
\hline
$E_7$ & $p = 19$ & $2p - 9 = 27$ \\
\hline
$E_8$ & $p \geq 31$ & $2p - 3$ \\
\hline
$F_4$ & $p \geq 13$ & $2p - 9 \leq D \leq 2p - 3$\\
\hline
$G_2$ & $p \geq 7$ & $2p - 8 \leq D \leq 2p - 3$\\
\hline
\end{tabular}
\end{center}

\vskip.2cm
For $r > 1$, one can usually obtain bounds inductively, with the sharp bound given by multiplying $D$ as above by $r$.  E.g., $r(2p-3)$ or $r(2p - 2n)$.  However, this is not possible in all cases, particularly types $B_n$ and $E_6$, where additional subtleties arise.  The reader is referred to \cite{BNP11, BNP12} for precise details on the known bounds.

The aforementioned work of Sprehn \cite{S} is purely constructive in nature, showing the existence of cohomology classes in certain degrees (but not the vanishing of any cohomology).  Sprehn considers an arbitrary connected reductive group $G$ and formally works over the finite field $\fp$ assuming that $p \geq h$.   In that work, it is shown that $\opH^{\bullet}(\gfq,\fp)$ contains (as a graded module over the Steenrod algebra) $\opH^{\bullet}(GL_2(\fq),\fp)$, from which it follows that $\opH^{r(2p-3)}(\gfq,\fp) \neq 0$.  By similar means, the lower degree non-vanishing is also observed for the symplectic group (type $C_n$): $\opH^{r(p-2)}(Sp_{2n}(\fq),\fp) \neq 0$ (for $p \geq 2n$).  

\subsection{Adjoint type}  Suppose for this subsection that $G$ is simple of {\em adjoint} type rather than simply connected, so that the root lattice coincides with the weight lattice. For example, we are considering $PGL_n(\fq)$ rather than $SL_n(\fq)$ in type $A$.  For the simply-laced root systems $A_n$, $D_n$, and $E_n$, still with $p > h$, there is a uniform sharp vanishing bound of $r(2p - 3)$ (cf. \cite{BNP12}).

\subsection{Small primes}\label{S:smallvan} In addition to the remaining questions for large primes, the task of identifying precise vanishing ranges for $p < h$ remains much more open.  Using a more direct approach to showing the vanishing of $\Ext^{j}_{G}(V(\la)^{(1)},H^0(\la))$, Bendel, Nakano, and Pillen \cite{BNP14} obtained smaller vanishing ranges for arbitrary primes:
\begin{itemize}
\item $p = 2$: $\opH^i(\gfq,k) = 0$ for $0 < i < r$,
\item $p > 2$: $\opH^i(\gfq,k) = 0$ for $0 < i < r(p-2)$.
\end{itemize}
For an arbitrary (simple, simply connected) $G$, those are the best possible ranges as seen by long-known computations of Carlson \cite{C} for $SL_2$: for $p = 2$, $\opH^r(SL_2(\fq),k) \neq 0$, and, for $p > 2$, $\opH^{r(p-2)}(SL_2(\fq),k) \neq 0$.   These ranges are of course not necessarily strict in general.

Following the characteristic class approach of \cite{S}, Lahtinen and Sprehn \cite{LS} gave a generalization in a sense of the aforementioned generic $r(2p-3)$ non-vanishing result. For $n \geq 2$, let $s$ be such that $p^{s-1} < n \leq p^{s}$.  They showed that $\opH^d(GL_n(\fq),\fp) \neq 0$ for $d = r(2p^s - 2p^{s-1} - 1)$.\footnote{They also gave the existence of non-nilpotent elements (in larger degree).}  Note that, when $s = 1$, $d = r(2p - 3)$.      For $r = 1$, some stronger results are given on the existence of non-nilpotent elements. For example, with $p = 2$, suppose $2 \leq n \leq 2^s$. If the sum of the binary digits of $d$ is at least $s$, then there is a non-nilpotent element in $\opH^d(GL_n(\fp),\fp)$.   For $p > 2$, the result involves the $p$-ary expansion of the degree.


\section{Finite group cohomology of simple modules}\label{S:small}

\subsection{} Similar to the situation with Frobenius kernels, computations of $\gfq$-cohomology with non-trivial coefficients are limited.   See \cite[\S 12 \& \S 14, particularly \S 12.11 \& \S 14.2]{Hum} for some early results.   We consider here the basic case of $\opH^i(\gfq,L(\la))$ for $\la \in X_r$.  Jones \cite{Jon}, Cline, Parshall, and Scott \cite{CPS75}, and Jones and Parshall \cite{JP} computed $\opH^1(\gfq,L(\la))$ for dominant weights $\la$ that are {\em minimal} in the usual partial ordering of weights.

\subsection{Degree 1} More recent work of the University of Georgia VIGRE Algebra Group \cite{UGA1} took a new approach to these degree 1 computations and expanded on the known results.  The approach of that work (related to the ideas in Section \ref{S:rat}) was to relate $\gfq$-cohomology to $G$-cohomology. Specifically they showed that for $q > 3$ and $p$ given in the following table

\vskip.2cm
\begin{center}
\begin{tabular}{|c|c|c|c|}
\hline
Root System & $A_n$, $D_n$ & $B_n$, $C_n$, $E_6$, $E_7$, $F_4$, $G_2$ & $E_8$\\
\hline
$p > $ & 2 & 3 & 5\\
\hline
\end{tabular},
\end{center}

\vskip.2cm\noindent
if $\la \leq \varpi_j$ for some fundamental dominant weight $\varpi_j$, then $\opH^i(\gfq,L(\la)) \cong \opH^1(G,L(\la))$.   This was done by relating $\gfq$-cohomology to $U_1$-cohomology.   With this result and a slightly stronger condition on the prime ($p > 7$ in types $E_7$ and $E_8$), the cohomology groups were computed and usually seen to be zero.\footnote{With exceptions in type $C_n$, $F_4$, $E_7$, and $E_8$, where the result is the trivial module $k$.}  In particular, under the assumptions on $q$ and $p$ as in the above table, complete computations of $\opH^1(\gfq,L(\varpi_j))$ are given.

\subsection{Degree 2} Using the induction functor $\ind_{\gfq}^G(-)$ and the ideas of Section \ref{S:genrev} and (again) connections with $U_r$-extensions, in later work, the University of Georgia VIGRE Algebra group \cite{UGA2} investigated degree 2 cohomology and identified conditions under which $\opH^2(\gfq,L(\la)) \cong \opH^2(G,L(\la))$.  Again, the focus is on the case that $\la \leq \varpi_j$. Such an isomorphism is given under the following assumptions on $p$ and $q$\footnote{$\dagger$ In type $C_n$, if $\la$ lies in the root lattice, one requires $q > 5$.}

\vskip.2cm
\begin{center}
\begin{tabular}{|c|c|c|c|}
\hline
Root System & $A_n$, $B_n$, $C_n^{\dagger}$, $D_n$, $E_6$ & $E_7$, $F_4$ & $G_2$\\
\hline
$p > $ & 3 & 3 & 5\\
\hline
$q > $ & 3  & 5 & 5\\
\hline
\end{tabular},
\end{center}

\vskip.2cm\noindent
with some exceptions on the weight $\la$.   Explicit computations are made in almost all cases, as well as computations for most of the exceptional cases. These latter cases give some examples where $\opH^2(G,L(\la)) = 0$ but $\opH^2(\gfq,L(\la)) \neq 0$.


\section{Bounding Dimensions}\label{S:bound}

\subsection{} For a finite simple group, a longstanding question has been to bound the size of the first cohomology group of a faithful, absolutely irreducible module.  In our context, the specific question is to understand the dimension of $\opH^1(\gfq,L(\si))$ for $\si \in X_1$.    Even in the broader finite simple group context, known examples (such as those in Section \ref{S:small}) suggested that the size should be small, potentially with a universal bound.  That remained the consensus until a workshop held at the American Institute of Mathematics in 2012, where discussions at the meeting and soon thereafter led to the discovery of extremely large-dimensional cohomology groups.   For a discussion of some of the history and development of the problem, the reader is referred to \cite{AIM} (cf. also \cite{BNPPSS}).  

The existence of these large-dimensional cohomology groups fundamentally changed the perspective on the problem, in particular, increasing the importance of earlier work of Cline, Parshall, and Scott \cite{CPS09} who had shown that there is a bound
$$
\dim\opH^1(\gfq,L(\la)) \leq C(\Phi),
$$
where $C(\Phi)$ depends only on the root system, as well as similar work of Guralnick and Tiep \cite{GT}, where the {\em cross-characteristic} case was considered. That is, where one is considering $\gfq$ with $p$ (the characteristic of the base field $k$) {\em not} dividing $q$; a situation generally not under consideration in this paper.  


\subsection{$G$-cohomology} The Cline-Parshall-Scott result along with further developments below provide another example of the key interplay between cohomology groups over different structures.  The bound on $\opH^1(\gfq,L(\la))$ was achieved by showing that such a bound held for $G$-cohomology of simple modules, and then using a result of Bendel, Nakano, and Pillen \cite{BNP06} (as referenced in Section \ref{S:shift}) that related $\gfq$-cohomology with $G$-cohomology.

Later work of Parshall and Scott \cite{PS11} extended these bounds on $G$-cohomology to higher degree cohomology groups, showing that there exist bounds on $\dim\opH^i(G,L(\la))$ depending only on $\Phi$ and the degree $i$.  More generally, it was shown that such a bound exists for $\dim\Ext_G^i(L(\la),L(\mu))$ for a pair of simple modules.

\subsection{Finite groups continued} Parker and Stewart \cite{PaS} improved upon \cite{CPS09} by finding explicit bounds on the dimension of $\opH^1(\gfq,L(\la))$ based on the Coxeter number $h$.  That work also makes use of the relationship with $G$-cohomology but takes a more direct approach to bounding the $G$-cohomology through the use of the Jantzen sum formula \cite[II.8.19]{Jan} for a Weyl module and bounding the number of composition factors in a Weyl module.  The explicit nature of their bound also allows one to see how the bound grows as the rank grows.

\subsection{Finite group cohomology in higher degrees} Using the bounds on $G$-cohomology in higher degrees from \cite{PS11}, Bendel, Nakano, Parshall, Pillen, Scott, and Stewart \cite{BNPPSS} showed that there exist bounds
$$
\dim\opH^i(\gfq,L(\la)) \leq C(\Phi, i),
$$
where $C(\Phi,i)$ depends only on $\Phi$ and the degree $i$.  More generally, such a bound was given for extensions between simple $\gfq$-modules.  This made use of a modification of ideas discussed in Sections \ref{S:rat} and \ref{S:van}. 

For a simple $G$-module $L(\si)$, as in Section \ref{S:indgfq}, we know $$\opH^i(\gfq,L(\si)) \cong \opH^i(G,\ind_{\gfq}^G(L(\si)) \cong \opH^i(G,L(\si)\otimes\ind_{\gfq}^G(k))$$ and $\ind_{\gfq}^G(k)$ admits a filtration with factors of the form $H^0(\la)\otimes H^0(\la^*)^{(r)}$, one for each $\la \in X^+$.    For a given $i$, there exists a positive integer $b$ such that $\opH^j(G,L(\si)\otimes H^0(\la)\otimes H^0(\la^*)^{(r)}) = 0$ for $j \leq i$ and $\langle \la, \al_0^{\vee}\rangle > b$. This allows one to truncate $\ind_{\gfq}^G(k)$ giving a finite-dimensional module $I_b$ with 
$$
\opH^i(\gfq,L(\si)) \cong \opH^i(G,L(\si)\otimes I_b),
$$
where $I_b$ admits a filtration with factors $H^0(\la)\otimes H^0(\la^*)^{(r)}$ where $\langle \la,\al_0^{\vee}\rangle \leq b$.  This then allows one to transfer $G$-bounds to $\gfq$.  

\begin{remark} In the \cite{BNPPSS} work, it was similarly shown that bounds exist on $\opH^i(G_r,L(\si))$ and, more generally, for Ext-groups between simple modules over Frobenius kernels.  
\end{remark}

\begin{remark} Given the potentially large dimensions of cohomology groups, it is of interest to determine just how large they can be.  Stewart \cite{Ste} demonstrated that the dimensions could be exponential in the cohomology degree by explicitly constructing a simple $SL_2$-module $L$ (depending on $n$) with $\dim\opH^{2n}(SL_2,L) \geq 2^{n-1}$.
\end{remark}

\subsection{Bounds based on the dimension of the coefficient module}  In a different direction from finding absolute bounds on dimensions, there have been attempts to bound cohomology in terms of the dimension of the coefficient module. That is, for a group $H$ and an $H$-module $M$, one may try to show that there is a constant $C$ (likely depending on the degree $i$) such that
$$
\dim\opH^i(H,M) \leq C\cdot\dim M.  
$$
For example, let $H$ be a finite group and $M$ be an irreducible $kH$-module on which $H$ acts faithfully, then Guralnick and Hoffman \cite{GH} showed that
$$
\dim\opH^1(H,M) \leq \frac12\cdot\dim M.  
$$
In degree 2, for a quasi-simple group $H$ and any finite-dimensional $kH$-module $M$, Guralnick, Kantor, Kassabov, and Lubotsky \cite{GKKL} showed that 
$$
\dim\opH^2(H,M) \leq 17.5\cdot\dim M.  
$$
For an arbitrary finite group and faithful, irreducible module $M$, they found a bound of $C = 18.5$.

In our setting, this question was investigated with some progress at the aforementioned 2012 AIM meeting by Bendel, Boe, Drupieski, Nakano, Parshall, Pillen, and Wright \cite{AIM}.   For the algebraic group $G$ and $M$ an arbitrary finite-dimensional rational $G$-module, the following bounds were found:
\begin{itemize}
\item $\dim\opH^1(G,M) \leq \frac{1}{h}\cdot\dim M \leq \frac12\cdot \dim M$,
\item $\dim\opH^2(G,M) \leq \dim M$,
\item for $p > h$, $\dim\opH^3(G,M) \leq 2\cdot\dim M$.
\end{itemize}
These results use the fact that $\opH^i(G,M) \cong \opH^i(B,M)$ for a rational $G$-module $M$ (as noted in Section \ref{S:Bintro}) and results from Section \ref{S:Bla} on $\opH^i(B,\la)$. In the degree 3 case, the condition $p > h$ can be weakened to those given in Section \ref{S:Lie} by using the $\opH^3(B,\la)$ computations of \cite{BNP16}.  Using these $G$-cohomology bounds, one may show that, by choosing $r$ sufficiently large,  the same bounds may be given for $\gfq$ (and any $k\gfq$-module $M$).  

\begin{remark} Outside the context of this paper, in \cite{AIM}, bounds were also obtained on the first cohomology of a symmetric group with coefficients in a simple module.  
\end{remark}

\providecommand{\bysame}{\leavevmode\hbox
to3em{\hrulefill}\thinspace}


\begin{thebibliography}{88888888888888}

\bibitem[\sf And]{And} H.~H.~Andersen, {Extensions of modules for algebraic groups}, {\em Amer. J. Math.}, {\sf 106}, (1984), 498--504.

\bibitem[\sf AJ]{AJ} H.~H.~Andersen, J.~C.~Jantzen, Cohomology of induced representations for algebraic groups, {\em Math. Ann.}, {\sf 269}, (1984), 487--525.

\bibitem[\sf AR]{AR} H.~H.~Andersen, T.~Rian, {$B$-cohomology}, {\em J. Pure Appl. Algebra}, {\sf 209}, (2007), 537--549.

\bibitem[\sf Bar]{B} A.~Barbu, On the range of non-vanishing $p$-torsion cohomology for $GL_n(\fp)$, {\em J. Algebra}, {\sf 278}, (2004), 456-472.

\bibitem[\sf AIM]{AIM} C.~P.~Bendel, B.~D.~Boe, C.~M.~Drupieski, D.~K.~Nakano, B.~J.~Parshall, C.~Pillen, C.~B.~Wright, Bounding the dimensions of rational cohomology groups, {\em Developments and retrospectives in Lie theory}, 51–69, Dev. Math., 38, Springer, Cham, 2014.


\bibitem[\sf BNPP]{BNPP} C.~P.~Bendel, D.~K.~Nakano, B.~J.~Parshall, C.~Pillen, Cohomology for quantum groups via the geometry of the nullcone, {\em Mem. Amer. Math. Soc.}, {\sf 229}, no. 1077,  (2014).

\bibitem[\sf BNPPSS]{BNPPSS} C.~P.~Bendel, D.~K.~Nakano, B.~J.~Parshall, C.~Pillen, L.~L.~Scott, D.~Stewart, Bounding cohomology for finite groups and Frobenius kernels, {\em Algebr. Represent. Theory}, {\sf 18}, (2015), 739–760. 

\bibitem[\sf BNP01]{BNP01} C.~P.~Bendel, D.~K.~Nakano, C.~Pillen, On comparing the cohomology of algebraic groups, finite Chevalley groups, and Frobenius kernels, {\em J. Pure \& Appl. Algebra}, {\sf 163}, (2001), 119-146.

\bibitem[\sf BNP02]{BNP02} C.~P.~Bendel, D.~K.~Nakano, C.~Pillen, Extensions for finite Chevalley groups II, {\em Trans. Amer. Math. Soc.}, {\sf 354}, (2002), 4421--4454.

\bibitem[\sf BNP04a]{BNP040} C.~P.~Bendel, D.~K.~Nakano, C.~Pillen, Extensions for finite Chevalley groups I, {\em Adv. Math.}, {\sf 183}, (2004), 380--408.


\bibitem[\sf BNP04b]{BNP041} C.~P.~Bendel, D.~K.~Nakano, C.~Pillen, Extensions for finite groups of Lie type: twisted groups, {\em Finite groups 2003}, 29-46, Walter de Gruyter GmbH \& Co. KG, Berlin, 2004.

\bibitem[\sf BNP04c]{BNP04}C.~P.~Bendel, D.~K.~Nakano, C.~Pillen, Extensions for Frobenius kernels, {\em J. Algebra}, {\sf 272}, (2004), 476-511.

\bibitem[\sf BNP06]{BNP06} C.~P. Bendel, D.~K. Nakano, C.~Pillen, Extensions for finite groups of Lie Type II: filtering the truncated induction functor, {\em Cont. Math.}, {\sf 413}, (2006), 1-23.

\bibitem[\sf BNP07]{BNP07} C.~P.~Bendel, D.~K.~Nakano, C.~Pillen, Second cohomology groups for Frobenius kernels and related structures, {\em Adv. Math.}, {\sf 209}, (2007), 162-197.

\bibitem[\sf BNP11]{BNP11} C.~P.~Bendel, D.~K.~Nakano, C.~Pillen, On the vanishing ranges for the cohomology of finite groups of Lie type, {\em Int. Math. Res. Not. IMRN 2012}, no. 12, 2817-2866; doi: 10.1093/imrn/rnr130.

\bibitem[\sf BNP12]{BNP12} C.~P.~Bendel, D.~K.~Nakano, C.~Pillen, On the vanishing ranges for the cohomology of finite groups of Lie type II, {\em  Recent developments in Lie algebras, groups and representation theory}, 25–73, Proc. Sympos. Pure Math., 86, Amer. Math. Soc., Providence, RI, 2012.

\bibitem[\sf BNP14]{BNP14} C.~P.~Bendel, D.~K.~Nakano, C.~Pillen, Extensions for finite Chevalley groups III: Rational and generic cohomology, {\em Adv. Math.}, {\sf 262}, (2014), 484–519.

\bibitem[\sf BNP16]{BNP16} C.~P.~Bendel, D.~K.~Nakano, C.~Pillen, Third cohomology groups for Frobenius kernels and related structures, {\em Lie algebras, Lie superalgebras, vertex algebras and related topics}, 81–118, Proc. Sympos. Pure Math., 92, Amer. Math. Soc., Providence, RI, 2016.

\bibitem[\sf BNPS20]{BNPS20} C.~P.~Bendel, D.~K.~Nakano, C.~Pillen, P.~Sobaje, Counterexamples to the tilting and $(p,r)$-filtration conjectures, {\em J. Reine Angew. Math}, {\sf 767}, (2020), 193-202.

\bibitem[\sf BNPS22a]{BNPS22} C.~P.~Bendel, D.~K.~Nakano, C.~Pillen, P.~Sobaje, On Donkin's tilting module conjecture II: Counterexamples, arXiv: 2107.11615. 

\bibitem[\sf BNPS22b]{BNPS23} C.~P.~Bendel, D.~K.~Nakano, C.~Pillen, P.~Sobaje, On Donkin's tilting module conjecture III: New generic lower bounds, arXiv: 2209.04675. 


\bibitem[\sf Car]{C} J.~F.~Carlson, The cohomology of irreducible modules over $SL(2, p^n)$, {\em Proc. London Math. Soc.} (3), {\sf 47}, (1983), 480–492.

\bibitem[\sf CPS75]{CPS75} E.~Cline, B.~ Parshall, L.~Scott, Cohomology of finite groups of Lie type. I,
{\em Inst. Hautes {\' E}tudes Sci. Publ. Math.}, {\sf 45}, (1975), 169–191.


\bibitem[\sf CPS09]{CPS09} E.~T.~Cline, B.~J.~ Parshall, L.~L.~Scott, Reduced standard modules and cohomology, {\em Trans. Amer. Math. Soc.}, {\sf 361}, (2009), 5223-5261.

\bibitem[\sf CPSvdK]{CPSvdK} E.~Cline, B.~Parshall, L.~Scott, W.~van der Kallen, Rational and generic cohomology, {\em Invent. math.}, {\sf 39}, (1977), 143–163.

\bibitem[\sf Don]{Don} S.~Donkin, On tilting modules for algebraic groups, {\em Math. Z.}, {\sf 212}, (1993), 39-60.


\bibitem[\sf DS]{DS} M.~F.~Dowd, P.~Sin, On representations of algebraic groups in characteristic two, {\em Comm. Algebra}, {\sf 24}, (1996), 2597–2686.


\bibitem[\sf Fr76]{F76} E.~M.~Friedlander, Computations of $K$-theories of finite fields, {\em Topology}, {\sf 15}, (1976), 87-109.


\bibitem[\sf Fr19]{F19} E.~M.~Friedlander, Cohomology of unipotent group schemes, {\em Algebr. Represent. Theory}, {\sf 22}, (2019), 1427-1455.

\bibitem[\sf FP83]{FP83} E.~M.~Friedlander, B.~J.~Parshall, On the cohomology of algebraic and related finite groups, {\em Invent. Math.}, {\sf 74}, (1983), 85-118.


\bibitem[\sf FP86a]{FP86} E.~M.~Friedlander, B.~J.~Parshall, Cohomology of Lie algebras and algebraic groups, {\em Amer. J. Math.}, {\sf 108}, (1986), 235-253.

\bibitem[\sf FP86b]{FP862} E.~M.~Friedlander, B.~J.~Parshall, Cohomology of infinitesimal and discrete groups, {\em Math. Ann.}, {\sf 273}, (1986), 353-374.


\bibitem[\sf FS]{FS} E.~M. Friedlander, A. Suslin, Cohomology of finite group schemes over a field, {\em Invent. Math.}, 127 (1997), 209-270.


\bibitem[\sf GH]{GH} R.~M. Guralnick, C.~F.~Hoffman, The first cohomology group and generation of simple groups, {\em Groups and geometries (Siena, 1996)}, 81-89, Trends Math., Birkh{\" a}user, Basel, 1998.

\bibitem[\sf GKKL]{GKKL} R.~M.~Guralnick, W.~M.~Kantor, M.~Kassabov, A.~Lubotzky, Presentations of finite simple groups: profinite and cohomological approaches, {\em Groups Geom. Dyn.}, {\sf 1}, (2007), 469-523.

\bibitem[\sf GT]{GT} R.~M.~Guralnick, P.~H.~Tiep, First cohomology groups of Chevalley groups in cross characteristic, {\em Ann. of Math.} (2), {\sf 174}, (2011), 543–559.



\bibitem[\sf Hil]{H} H.~L.~Hiller, Cohomology of Chevalley groups over finite fields, {\em J. Pure Appl. Algebra}, {\sf 16}, (1980), 259-263.


\bibitem[\sf Hum]{Hum} J.~E.~Humphreys, {\em Modular Representations of Finite Groups of Lie Type}, London Mathematical Society Lecture Notes Series, {\sf 326}, Cambridge University Press, Cambridge, 2006.

\bibitem[\sf HV]{HV} J.~E.~Humphreys, D.~N.~Verma, Projective modules for finite Chevalley groups, {\em Bull. Amer. Math. Soc.}, {\sf 79}, (1973), 467-468.

\bibitem[\sf Jan91]{Jan91} J.~C.~Jantzen, First cohomology groups for classical Lie algebras, {\em Representation Theory of Finite Groups and Finite-Dimensional Algebras (Bielefeld, 1991)}, 289-315, Progr. Math., {\sf 95}, Birkh{\" a}user, Basel, 1991.

\bibitem[\sf Jan03]{Jan} J.~C.~Jantzen, {\em Representations of Algebraic Groups}, Mathematical Surveys and Monographs, {\sf 107}, American Mathematical Society, Providence, RI, 2003.

\bibitem[\sf Jon]{Jon} W.~Jones, {\em Cohomology of finite groups of Lie type}, Ph.D. thesis, University of Minnesota,
1975.

\bibitem[\sf JP]{JP} W.~Jones, B.~Parshall, On the 1-cohomology of finite groups of Lie type, {\em Proceedings
of the Conference on Finite Groups (Univ. Utah, Park City, Utah, 1975)}, 313-328,  Academic Press, New York, 1976. 

\bibitem[\sf KSTY]{KSTY} M.~Kaneda, N.~Shimada, M.~Tezuka, N.~Yagita, Cohomology of infinitesimal algebraic groups, {\em Math. Z.}, {\sf 205}, (1990), 61–96.


\bibitem[\sf Kos]{Kos} B.~Kostant, Lie algebra cohomology and the generalized Borel-Weil theorem, {\em Ann. Math.}, {\sf 74}, (1961), 329-387.

\bibitem[\sf KLT]{KLT} S.~Kumar, N.~Lauritzen, J.~Thomsen, Frobenius splitting of cotangent bundles of flag varieties, {\em Invent. Math.}, {\sf 136}, (1999), 603-621.


\bibitem[\sf LS]{LS} A.~Lahtinen, D.~Sprehn, Modular characteristic classes for representations over finite fields, {\em Adv. Math.}, {\sf 323}, (2018), 1–37.

\bibitem[\sf Lin]{Lin} Z.~Lin, Extensions between simple modules for Frobenius kernels, {\em Math. Z.}, {\sf 207}, (1991), 485-499.

\bibitem[\sf Nak]{Nak} D.~K.~Nakano, Cohomology of algebraic groups, finite groups, and Lie algebras: interactions and connections, {\em Lie theory and representation theory}, 151–176, Surv. Mod. Math., {\sf 2}, Int. Press, Somerville, MA, 2012.


\bibitem[\sf Ngo13]{N13} N.~Ngo,  Cohomology for Frobenius kernels of $SL_2$, {\em J. Algebra}, {\sf 396}, (2013), 39-60.

\bibitem[\sf Ngo19]{N19} N.~Ngo, Low degree cohomology of Frobenius kernels, {\em Advances in algebra}, 255–259,
Springer Proc. Math. Stat., 277, Springer, Cham, 2019.

\bibitem[\sf OHal]{OHal} J.~O'Halloran, A vanishing theorem for
the cohomology of Borel subgroups, {\em Comm. Algebra}, {\sf 11}, (1983),
1603-1606.

\bibitem[\sf PaSt]{PaS} A.~Parker, D.~I.~Stewart, First cohomology groups for finite groups of Lie type in defining characteristic, {\em Bull. London Math. Soc.}, {\sf 46}, (2014), 227–238.

\bibitem[\sf PS]{PS11} B.~J.~Parshall, L.~L.~Scott, Bounding Ext for modules for algebraic groups, finite groups and quantum groups, {\em Adv. Math.}, {\sf 226}, (2011), 2065–2088.

\bibitem[\sf PSS]{PSS} B.~J.~Parshall, L.~L.~Scott, D.~I.~Stewart, Shifted generic cohomology, {\em Compos. Math.}, {\sf 149}, (2013), 1765-1788.

\bibitem[\sf PT]{PT} P.~Polo, J.~Tilouine, {Berstein-Gelfand-Gelfand complexes and cohomology of nilpotent groups over $\mathbb{Z}_{(p)}$ for representations with $p$-small weights}, {\em Ast{\'e}risque}, {\sf 280}, (2002), 97-135.

\bibitem[\sf Rad]{Rad} A.-C.~Radu, On the cohomology of the Ree groups and kernels of exceptional isogenies, {\em J. Algebra}, {\sf 612}, (2022), 1-54.

\bibitem[\sf RW]{RW} S.~Riche, G.~Williamson, A simple character formula, {\em Ann. H. Lebesgue}, {\sf 4}, (2021), 503–535.

\bibitem[\sf Qui]{Q} D.~G.~Quillen, On the cohomology and $K$-theory of finite fields, {\em Ann. Math.}, {\sf 96}, (1972), 551-586.

\bibitem[\sf Sin]{Sin} Peter Sin, Extensions of Simple modules for special algebraic groups, {\em J. Algebra}, 170, (1994), 1011-1034.

\bibitem[\sf Spr]{S} D.~Sprehn, Nonvanishing cohomology classes on finite groups of Lie type with Coxeter number at most $p$, {\em J. Pure Appl. Algebra}, {\sf 219}, (2015), 2396–2404.

\bibitem[\sf Ste12a]{Ste} D.~I.~Stewart, Unbounding Ext, {\em J. Algebra}, {\sf 365}, (2012), 1–11.

\bibitem[\sf Ste12b]{Ste2} D.~I.~Stewart, The second cohomology of simple $SL_3$-modules, {\em Comm. Algebra}, {\sf 40}, (2012), 4702–4716. 



\bibitem[\sf SFB1]{SFB1} A.~Suslin, E.~M.~Friedlander, C.~P.~Bendel,
Infinitesimal 1-parameter subgroups and cohomology, {\em Jour. Amer.
Math. Soc.}, {\sf 10}, (1997), 693-728.


\bibitem[\sf SFB2]{SFB2}  A.~Suslin, E.~M.~Friedlander, C.~P.~Bendel,
Support varieties for infinitesimal group schemes, {\em Jour. Amer.
Math. Soc.}, {\sf 10}, (1997), 729-759.

\bibitem[\sf UGA09]{UGA} University of Georgia VIGRE Algebra Group,
On Kostant's theorem for Lie algebra cohomology, {\em Cont. Math.}, 
{\sf 478}, (2009), 39--60. 

\bibitem[\sf UGA12]{UGA2} University of Georgia VIGRE Algebra Group,
Second cohomology for finite groups of Lie type,
{\em J. Algebra}, {\sf 360}, (2012), 21-52.

\bibitem[\sf UGA13]{UGA1} University of Georgia VIGRE Algebra Group,
First cohomology for finite groups of Lie type: simple modules with small dominant weights,
{\em Trans. Amer. Math. Soc.}, {\sf 365}, (2013), 1025-1050.


\bibitem[\sf Wri]{W} C.~B.~Wright, Second cohomology groups for algebraic groups and their Frobenius kernels, {\em J. Algebra}, {\sf 330}, (2011), 60-75.

\bibitem[\sf Yeh]{Yeh} Samy El Badawy Yehia, {\em Extensions of simple modules for the universal Chevalley group and its parabolic subgroups}, Ph.D. thesis, University of Warwick, 1982.


\end{thebibliography}
\end{document}